\documentclass[12pt]{article}
\usepackage[english]{babel}
\usepackage{lineno} 
\usepackage{graphicx,multicol}
\usepackage{epic,eepic,epsfig}
\usepackage{caption}

\usepackage{amsmath,amsfonts,amsthm,amssymb}
\usepackage{pstricks,pstricks-add,pst-math,pst-xkey}

\usepackage{algorithm, algorithmic}

\usepackage{graphicx,pstricks,pst-node,amssymb}

\usepackage{tikz}

\pgfdeclarelayer{edgelayer}
\pgfdeclarelayer{nodelayer}
\pgfsetlayers{edgelayer,nodelayer,main}

\tikzstyle{none}=[inner sep=0pt]
\definecolor{hexcolor0xf81e1c}{rgb}{0.973,0.118,0.110}
\definecolor{hexcolor0x3c00ff}{rgb}{0.235,0.000,1.000}

\tikzstyle{whitevertex}=[circle,fill=white,draw=black, scale = 0.5]
\tikzstyle{redvertex}=[circle,fill=hexcolor0xf81e1c,draw=black, scale = 0.5]
\tikzstyle{bluevertex}=[circle,fill=hexcolor0x3c00ff,draw=black, scale = 0.5]
\tikzstyle{greenvertex}=[circle,fill=green,draw=black, scale=0.5]
\tikzstyle{purplevertex}=[circle,fill=magenta,draw=black, scale=0.5]
\tikzstyle{grayvertex}=[circle,fill=white,draw=gray, scale=0.5]
\tikzstyle{blackvertex}=[circle,fill=black,draw=black, scale=0.5]

\tikzstyle{textbox}=[rectangle,fill=none,draw=none]
\tikzstyle{box}=[rectangle,fill=none,draw=black]

\tikzstyle{arc}=[black, ->]
\tikzstyle{grayarc}=[gray, ->]
\tikzstyle{bluearc}=[blue, ->]
\tikzstyle{grayedge}=[draw=gray]
\tikzstyle{blueedge}=[draw=blue]
\tikzstyle{rededge}=[draw=red]
\tikzstyle{edge}=[draw=black]

\tikzstyle{vertex}=[circle, ,fill=white,draw=black, scale=0.5]

\tikzstyle{10circle}=[circle, scale=10.0,draw=black]
\tikzstyle{10oval}=[ellipse, scale=10.0,draw=black]

\usepackage{amsmath,amsfonts,amsthm,amssymb}
\usepackage{pstricks,pstricks-add,pst-math,pst-xkey}

\usepackage{algorithm, algorithmic}

\usepackage{graphicx,pstricks,pst-node,amssymb}

\usepackage{algorithm, algorithmic}

\begin{document}
\date{\today}

\newtheorem{tm}{\hspace{5mm}Theorem}
\newtheorem{prp}[tm]{\hspace{5mm}Proposition}
\newtheorem{dfn}[tm]{\hspace{5mm}Definition}
\newtheorem{lemma}[tm]{\hspace{5mm}Lemma}
\newtheorem{cor}[tm]{\hspace{5mm}Corollary}
\newtheorem{conj}[tm]{\hspace{5mm}Conjecture}
\newtheorem{prob}[tm]{\hspace{5mm}Problem}
\newtheorem{quest}[tm]{\hspace{5mm}Question}
\newtheorem{alg}[tm]{\hspace{5mm}Algorithm}
\newtheorem{sub}[tm]{\hspace{5mm}Algorithm}
\newcommand{\induce}[2]{\mbox{$ #1 \langle #2 \rangle$}}
\newcommand{\2}{\vspace{2mm}}
\newcommand{\dom}{\mbox{$\rightarrow$}}
\newcommand{\ndom}{\mbox{$\not\rightarrow$}}
\newcommand{\compdom}{\mbox{$\Rightarrow$}}
\newcommand{\cdom}{\compdom}
\newcommand{\sdom}{\mbox{$\Rightarrow$}}
\newcommand{\lsd}{locally semicomplete digraph}
\newcommand{\lt}{local tournament}
\newcommand{\la}{\langle}
\newcommand{\ra}{\rangle}
\newcommand{\pf}{{\bf Proof: }}
\newtheorem{claim}{Claim}
\newcommand{\beq}{\begin{equation}}
\newcommand{\eeq}{\end{equation}}
\newcommand{\<}[1]{\left\langle{#1}\right\rangle}

\newcommand{\Z}{\mathbb{$Z$}}
\newcommand{\Q}{\mathbb{$Q$}}
\newcommand{\R}{\mathbb{$R$}}


\title{Strong cocomparability graphs and Slash-free orderings of matrices}

\author{Pavol Hell\thanks{School of Computing Science, Simon Fraser University,
Burnaby, B.C., Canada V5A 1S6; pavol@sfu.ca}, 
Jing Huang\thanks{Department of Mathematics and Statistics,
      University of Victoria, Victoria, B.C., Canada V8W 2Y2; huangj@uvic.ca} and 
Jephian C.-H. Lin\thanks{Department of Applied Mathematics,
      National Sun Yat-sen University, Kaohsiung 80424, Taiwan; 
      jephianlin@gmail.com}}
\date{}

\maketitle

\begin{abstract}
We introduce the class of strong cocomparability graphs, as the class of reflexive 
graphs whose adjacency matrix can be rearranged by a simultaneous row and column 
permutation to avoid the submatrix with rows $01,10$, which we call {\em Slash}.

We provide an ordering characterization, a forbidden structure characterization, 
and a polynomial-time recognition algorithm, for the class. These results complete 
the picture in which in addition to, or instead of, the $Slash$ matrix one forbids the 
{\em $\Gamma$} matrix (which has rows $11,10$). It is well known that in these 
two cases one obtains the class of interval graphs, and the class of strongly 
chordal graphs, respectively.

By complementation, we obtain the class of strong comparability graphs, whose
adjacency matrix can be rearranged by a simultaneous row and column permutation 
to avoid the two-by-two identity submatrix. Thus our results give characterizations
and algorithms for this class of irreflexive graphs as well. In other words, our results
may be interpreted as solving the following problem: given a symmetric $0,1$-matrix
with $0$-diagonal, can the rows and columns of be simultaneously permuted to avoid 
the two-by-two identity submatrix?
\end{abstract}


\section{Introduction}

The problem of deciding whether a given $0,1$-matrix $M$ can have its rows and 
columns independently permuted to avoid containing a fixed submatrix $F$ (or a 
fixed family $\mathcal F$ of such submatrices) has been much studied in 
the literature, see for instance 
\cite{bd,fhh,gg,hhlm,mastro,huang,huang-oo,hks,krw,lubiw1,lub,pt,sbs}. 

Specifically, we say that a matrix $M$ is $F${\em-free}, or {\em avoids $F$}, if it
does not include a copy of $F$ as a submatrix. An $F${\em-free ordering} of a matrix
$M$ is an ordering of its rows and an ordering of its columns such that the 
resulting permuted matrix is $F$-free. The research mentioned above is concerned 
with deciding whether a given $M$ admits an $F$-free ordering. Often these problems
have interesting interpretations for graphs, and especially for bipartite graphs 
(i.e., {\em bigraphs}). Indeed, any $0,1$-matrix is associated with a bigraph as 
follows. The {\em biadjacency matrix} of a bigraph $G$ with parts $I, J$ is a 
$0,1,$-matrix $M$ whose rows are indexed by vertices of $I$ and columns by vertices
of $J$ and $M(i,j)=1$ if and only if $i \in I$ is adjacent to $j \in J$.

In the most widely known special case the matrix $F$ is the so-called 
{\em $\Gamma$ matrix}
$\left[\begin{matrix}1&1\\
                1&0
       \end{matrix}\right]$. 
It is known that a matrix $M$ admits a $\Gamma$-free ordering if and only if it is 
totally balanced (i.e., does not contain a cycle submatrix) \cite{af,hks,lub}),
and efficient algorithms are known to determine whether or not this is the case 
\cite{lub,pt}. Totally balanced matrices have important applications in location 
theory \cite{hks}. They also have a natural graph interpretation. We say that 
a bigraph is {\em chordal} if it does not contain an induced cycle of length 
at least six.
                
\begin{tm} \cite{af,hsk,lub}
A matrix admits a $\Gamma$-free ordering if and only if it is the biadjacency 
matrix of a chordal bigraph.
\qed
\end{tm}

This gives a forbidden structure characterization (in terms of the induced cycles 
of length at least six in the associated bigraph) of matrices that admit 
a $\Gamma$-free ordering.
                
We next review the situation when $F$ is the so-called {\em Slash} matrix 
$\left[\begin{matrix}0&1\\
                1&0
                \end{matrix}\right]$.
                
This matrix turns out to be important from the graph theoretic point of view, but 
we also point out that exchanging $0$'s and $1$'s translates this problem into one 
seeking $I_2$-free orderings, where $I_2$ is the two-by-two identity matrix. This
problem has attracted attention in the Russian literature because it yields 
efficiently solvable cases of a certain plant location problem \cite{krw}. It is 
known that $M$ admits an $I_2$-free ordering if and only if it is connected (i.e., 
the sequence of differences between the corresponding entries of any two rows 
changes sign at most once), and efficient recognition algorithms have been reported
in the Russian sources \cite{bd}.             

For the $Slash$ matrix, we have the following graph-theoretic interpretation. 
An {\em edge-asteroid} is a set of edges $x_0y_0, x_1y_1, \dots, x_{2k}y_{2k}$ such
that for each $i = 0, 1, \dots, 2k$, there is a walk that begins with the edge 
$x_{i+k}y_{i+k}$ and ends with the edge $x_{i+k+1}y_{i+k+1}$ such that neither 
$x_{i}$ nor $y_{i}$ is adjacent to a vertex in the walk. We call a bigraph 
a {\em cocomparability} bigraph if it contains no edge-asteroid \cite{hhlm}. 
This is motivated by an analogous definition and result of Gallai \cite{gallai}, 
stating that a graph is a comparability graph -- i.e., admits a transitive 
orientation -- if and only if its complement contains no asteroid. Here an asteroid
is a vertex-focused analogue of an edge-asteroid.; see \cite{hhlm} for details.

\begin{tm} \cite{hhlm} \label{hhlm}
A matrix admits a $Slash$-free ordering if and only if it is the biadjacency matrix
of a cocomparability bigraph.
\qed
\end{tm}

This result also gives a forbidden structure characterization of matrices that 
admit a $Slash$-free ordering (in terms of the edge-asteroids of the associated 
bigraph). It implies a robust polynomial-time recognition algorithm, which can
also be applied to recognizing matrices that admit $I_2$-free orderings.
The need for such algorithms is discussed in \cite{spin}, Section 10.3.

Finally, we review the situation for the family $\mathcal F = \{\Gamma,Slash\}$,
i.e., when both $\Gamma$ and $Slash$ are forbidden. (When the order of rows is
reversed, the resulting two matrices form a subset of $0,1$-Monge matrices, as
discussed in \cite{krw}, Sections 2.7 and 4.3.) There are several equivalent
descriptions of the class of matrices admitting $\mathcal F$-free orderings, 
in terms of bigraphs: as complements or circular arc graphs \cite{fhh}, as 
two-dimensional orthogonal ray graphs \cite{stu}, or as interval containment 
bigraphs \cite{hhmr,huang,huang-oo}, cf. also \cite{mastro}. We may state 
the theorem using already defined terms as follows \cite{fhh,huang}.

\begin{tm} \cite{hhmr,huang-oo} 
A matrix admits a $(\Gamma,Slash)$-free ordering if and only if it is the 
biadjacency matrix of a chordal bigraph without edge-asteroids.
\qed
\end{tm}

Other results give polynomial-time algorithms and characterizations in terms of 
bigraphs when similar sets of $2 \times 2$ matrices are forbidden \cite{krw,sbs}. 

All the results summarized above, translate to facts about classes of bigraphs. 
There is a more natural matrix candidate as far as correspondence to graphs is 
concerned, namely instead of the biadjacency matrix one can consider the usual 
adjacency matrix of an arbitrary graph. This is the focus of our paper. 

For technical reasons, we need to deal with {\em reflexive graphs} (graphs in which
ever vertex has a loop). In matrix terms, this means we have a symmetric 
$0,1$-matrix with $1$'s on the main diagonal. In the literature this is often 
presented without considering loops, but instead of the usual adjacency matrix one 
takes the {\em augmented adjacency matrix}, which is just the usual adjacency matrix
augmented with $1$'s on the main diagonal. Note that re-ordering the vertices of 
a graph corresponds to a simultaneous row and column permutation of its adjacency 
matrix.

It is also possible, of course, to consider an arbitrary square matrix to be 
the adjacency matrix of a digraph, and ask similar questions about avoiding given 
submatrices for simultaneous row and column permutations, see for instance 
\cite{hchl}. However, in this paper we focus purely on undirected graphs.

Thus from now on $M$ is a square symmetric $0,1$-matrix with $1$'s on the main 
diagonal. We now examine the situation for the same matrices $\Gamma$ and $Slash$.

A {\em symmetric ordering} of a symmetric matrix $M$ is a matrix obtained from $M$ 
by simultaneously permuting the rows and columns of $M$. Clearly, a symmetric 
ordering of a symmetric matrix is again symmetric. 

We first describe the situation when the $\Gamma$ matrix is forbidden.

A {\em trampoline} is a complete graph on $k$ vertices $u_1, u_2, \dots, u_k$ with 
$k \geq 3$, together with an independent set of $k$ vertices $v_1, v_2, \dots, v_k$
such that each $v_i$ is adjacent to $u_i$ and $u_{i+1}$ (and to no other vertices).
A graph is {\em strongly chordal} if it does not contain, as an induced subgraph, 
a cycle of length $\geq 4$ or a trampoline \cite{far}.

\begin{tm} \cite{far} \label{gamma}
A matrix admits a symmetric $\Gamma$-free ordering if and only if it is 
the adjacency matrix of a strongly chordal graph.
\qed
\end{tm}

This result provides a graph description of the class of matrices, but also 
a forbidden structure characterization, and a polynomial-time recognition algorithm.

For the family $\mathcal F = \{\Gamma,Slash\}$, the symmetric matrices that admit a 
$\mathcal F$-free ordering also have an elegant graph description, which results in
a well-known forbidden structure characterization and polynomial-time recognition.

A graph $G$ is an {\em interval} graph if there is a family of intervals 
$I_v, v \in V(G)$ such that for all $u, v \in V(G)$,  $uv \in E(G)$ if and only 
if $I_u \cap I _v \neq \emptyset$.  

\begin{tm} \cite{adjusted,hhmr} \label{gamma-slash}
A matrix admits a symmetric $(\Gamma,Slash)$-free ordering if and only if it is
the adjacency matrix of an interval graph.
\qed
\end{tm}

The purpose of our paper is to complete the picture by introducing the ``last piece 
of the puzzle'', that is, the class of reflexive graphs whose adjacency matrices 
have symmetric $Slash$-free orderings. We call such graphs 
{\em strong cocomparability graphs}. We provide a natural ordering characterization
of strong cocomparability graphs, as well as a characterization by forbidden 
substructures that is analogous to similar characterizations of interval graphs and
other graph and digraph classes. This latter characterization implies a 
polynomial-time algorithm for the recognition of strong cocomparability graphs.

We apply our results to prove that the class of interval graphs is the intersection
of the classes of strongly chordal graphs and strong cocomparability graphs, 
enhancing a well known theorem of Gilmore and Hoffman \cite{gh}. 

These results can be applied by complementation to the class of 
{\em strong comparability graphs}, which correspond to adjacency matrices that 
have a symmetric $I_2$-free ordering. We ask if the domination problem, NP-complete
for comparability graphs, can be solved in polynomial time for strong comparability
graphs, in analogy to the situation for chordal versus strongly chordal graphs.

Unfortunately, we do not have a forbidden subgraph characterization of strong 
cocomparability graphs, although such characterizations are known for interval
graphs, cocomparability graphs, chordal graphs, and strongly chordal graphs. 
We do, however, provide a simple description of the class of bipartite strong 
cocomparability graphs, and relate it to the class of cocomparability bigraphs and 
to the class of bipartite permutation graphs.

\section{Symmetric $Slash$-free orderings}

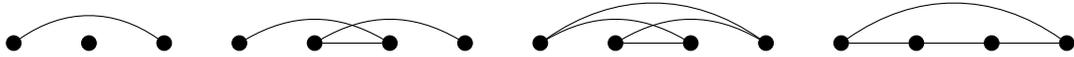
\begin{figure}[ht]
\begin{center}
\begin{tikzpicture}[>=latex]
                \node [style=blackvertex] (2) at (-5,0) {};
                \node [style=blackvertex] (3) at (-4,0) {};
                \node [style=blackvertex] (4) at (-3,0) {};
                \draw [ ] (2) to [out=35,in=145] (4);
\node [style=blackvertex] (5) at (-2,0) {};
                \node [style=blackvertex] (6) at (-1,0) {};
                \node [style=blackvertex] (7) at (0,0) {};
                \node [style=blackvertex] (8) at (1,0) {};
                \draw [ ] (5) to [out=30,in=150] (7);
                \draw [ ] (6) to [out=30,in=150] (8);
                \draw [ ] (6) to (7);
\node [style=blackvertex] (9) at (2,0) {};
                \node [style=blackvertex] (10) at (3,0) {};
                \node [style=blackvertex] (11) at (4,0) {};
                \node [style=blackvertex] (12) at (5,0) {};
                \draw [ ] (9) to [out=30,in=150] (11);
                \draw[ ] (10) to [out=30,in=150] (12);
  \draw[ ] (9) to [out=35,in=145] (12);
     \draw [ ] (11) to (10);
\node [style=blackvertex] (13) at (6,0) {};
                \node [style=blackvertex] (14) at (7,0) {};
                \node [style=blackvertex] (15) at (8,0) {};
                \node [style=blackvertex] (16) at (9,0) {};
                \draw [ ] (13) to [out=35,in=145] (16);
                \draw[ ] (13) to (14);
                \draw[ ] (14) to (15);
                \draw[ ] (16) to (15);

\end{tikzpicture}
\end{center}
\vspace{-2mm}
\caption{Forbidden patterns for strong cocomparability orderings.}
\label{forbpat}
\end{figure}

Let $\prec$ be a vertex ordering of a reflexive graph $G$ and $M$ be the adjacency
matrix of $G$ whose rows and columns are ordered according to $\prec$. Then
$M$ contains the $Slash$ matrix as a submatrix if and only if $\prec$ contains, as
ordered induced subgraphs, a pattern in Figure \ref{forbpat}. 
This can be easily seen by considering the possible relative positions of the row 
numbers and the column numbers of the occurrence of $Slash$.

Thus we have the following ordering characterization of strong cocomparability
graphs.

\begin{prp} \label{scoco-order}
A reflexive graph is a strong cocomparability graph if and only if it admits
a vertex ordering that does not contain any of the four patterns in 
Figure \ref{forbpat}.
\qed
\end{prp}

We call such a vertex ordering of $G$ a {\em strong cocomparability ordering}.
It follows by comparison with the cocomparability ordering mentioned in \cite{dama}
that a strong cocomparability graph is also a cocomparability graph.

We now proceed to give a structural characterization of strong cocomparability
graphs. Our result parallels other structural characterizations of various classes
of graphs and digraphs \cite{signed,adjusted,cags,mastro,boj,mono}, as it is given
in terms of absence of so-called invertible pairs, defined below. We point out that,
for example, interval graphs have a very similar structural characterization
\cite{adjusted}. To prove that the absence of invertible pairs implies the
existence of a strong cocomparability ordering, we use the method of
``complementary circuits" pioneered in \cite{adjusted}, and also used in
\cite{signed,boj,mono}.

In particular, our characterization will imply a polynomial-time recognition
algorithm.

We first introduce a useful ``forcing" relation $\Lambda$.

Let $G$ be a reflexive graph with vertex $V(G)$ and edge set $E(G)$. Note that 
$E(G)$ includes all {\em loops} $vv$, $v \in V(G)$. Denote by $Z(G)$ 
the set of ordered pairs $(u,v)$ of distinct vertices of $G$. 
For $(u,v), (u',v') \in Z(G)$, we say that $(u,v)$ {\em forces} $(u',v')$, 
denoted by $(u,v) \Lambda (u',v')$, if $u = u'$ and $v = v'$ or $uu', vv' \in E(G)$ 
and $uv', vu' \notin E(G)$. Clearly, $(u,v) \Lambda (u',v')$ if and only if 
$(v,u) \Lambda (v',u')$. We also say that $(u,v)$ {\em implies} $(u',v')$, denoted 
by $(u,v) \sim (u',v')$, if there exist walks
$u_1u_2 \dots u_k$ and $v_1v_2 \dots v_k$ in $G$ where
$(u_1,v_1) = (u,v)$ and $(u_k,v_k) = (u',v')$ such that
$(u_i,v_i) \Lambda (u_{i+1},v_{i+1})$ for each $i = 1, 2, \dots, k-1$. 
It is easy to verify that $\sim$ is an equivalence relation on $Z(G)$.
An equivalence class of $\sim$ will be called an {\em implication class}.
We shall use $[(u,v)]$ to denote the implication class of $(u,v)$ and
call $[(u,v)]$ {\em non-trivial} if $(u,v)$ is not the only pair in $[(u,v)]$.

The next proposition follows immediately from the definitions of 
strong cocomparability orderings and the relation $\Lambda$. 

\begin{prp} \label{pc}
Let $G$ be a reflexive graph. Then the following statements hold.
\begin{enumerate}
\item If $abc$ is an induced $P_3$ in $\overline{G}$, then $(a,b) \Lambda (c,b)$.
\item If $abcd$ is an induced $P_4$ in $G$, then
           $(a,c) \Lambda (a,d) \Lambda (b,d) \Lambda (a,d) \Lambda (b,c)$.
\item If $abcd$ is an induced $C_4$ in $G$, then $(a,d) \Lambda(b,c)$.
\qed
\end{enumerate}
\end{prp}

Each vertex ordering $\prec$ of $G$ gives rise to the set $\{(u,v):\ u \prec v\}$
which we denote by $U_{\prec}$. The set $U_{\prec}$ contains either 
$(u,v)$ or $(v,u)$ for all distinct vertices $u, v$ of $G$. Suppose that $\prec$ is 
a strong cocomparability ordering. Then for all $(u,v) \Lambda (u',v')$, 
$U_{\prec}$ either contains both $(u,v)$ and $(u',v')$ or both $(v,u)$ and 
$(v',u')$. Hence 
$U_{\prec}$ contains either $[(u,v)]$ or $[(v,u)]$ for all distinct vertices 
$u, v$ of $G$ and in particular, it is a union of implication classes. 
Conversely, if $U_{\prec}$ is a union of implication classes then $\prec$ is a 
strong cocomparability ordering of $G$.

\begin{prp}\label{scc}
Let $\prec$ be a vertex ordering of a reflexive graph $G$ and let
$U_{\prec} = \{(u,v) \in Z(G):\ u \prec v\}$.
The following statements are equivalent:
\begin{enumerate}
\item $\prec$ is a strong cocomparability ordering;
\item for all $(u,v) \sim (u',v')$, $u \prec v$ if and only if $u' \prec v'$;
\item $U_{\prec}$ is a union of implication classes.
\qed
\end{enumerate}
\end{prp}

An {\em invertible pair} in $G$ is a pair of distinct vertices $u, v$ such that
$(u,v) \sim (v,u)$. Clearly, if $u, v$ are an invertible pair then 
$[(u,v)] = [(v,u)]$ and for any $(x,y) \in [(u,v)]$, $x,y$ are an invertible pair.
Invertible pairs are absent in a strong cocomparability graph.

\begin{prp} \label{ip}
A strong cocomparability graph cannot contain an invertible pair.
\end{prp}

\pf If $G$ is a strong cocomparability graph which contains an invertible pair 
$u, v$, then by Proposition \ref{scc} any strong cocomparability ordering 
$\prec$ of $G$ would have $u \prec v$ and $v \prec u$, a contradiction.
\qed

An invertible pair in a reflexive graph $G$ (if exists) can be easily discovered 
via an auxiliary graph of $G$ defined next.

Define the {\em auxiliary graph} $G^+$ of reflexive graph $G$ as follows. 
The vertex set of $G^+$ is the set $Z(G)$.
Two vertices $(u,v)$ and $(v',u')$ are adjacent in $G^+$ if and only if
$(u,v) \Lambda (u',v')$, that is, $u = u'$ and $v = v'$, or $uu', vv' \in E(G)$ and 
$uv', vu' \notin E(G)$.

\begin{prp} \label{aux}
A reflexive graph $G$ has an invertible pair if and only if its auxiliary graph 
$G^+$ is not bipartite.
\end{prp}
\pf Suppose that $u_1, v_1$ are an invertible pair in $G$. Then there exist
walks $u_1u_2 \dots u_{k-1}u_k$ and $v_1v_2 \dots v_{k-1}v_k$ where 
$u_k = v_1$ and $v_k = u_1$ such that $(u_i,v_i) \Lambda (u_{i+1},v_{i+1})$ for
each $i = 1, 2, \dots, k-1$. Hence 
\[(u_1,v_1)(v_1,u_1)(u_2,v_2)(v_2,u_2) \cdots
(u_{k-1},v_{k-1})(v_{k-1},u_{k-1})(u_k,v_k)\]
is a closed walk of odd length in $G^+$, showing that $G^+$ is not bipartite.

Conversely, suppose that $G^+$ is not bipartite.
Let $(u_1,v_1)(u_2,v_2)\cdots (u_{2k+1},v_{2k+1})$ be a cycle of odd length
in $G^+$. Then 
\[(u_1,v_1) \Lambda (u_3,v_3) \Lambda \cdots \Lambda (u_{2k+1},v_{2k+1}) \Lambda 
(v_1,u_1)\]
showing that $u_1, v_1$ are an invertible pair in $G$.
\qed

In the next section we prove that the absence of invertible pairs in a reflexive 
graph in fact ensures that the graph is a strong cocomparability graph.

\section{The main result}

Let $G$ be a reflexive graph. For each implication class $I$ of $G$, denote by 
$I^{-1}$ the {\em inverse} class of $I$, that is,
\[I^{-1} = \{(u,v):\ (v,u) \in I\}. \]
A {\em circuit} of length $\ell$ is a sequence of pairs 
\[(u_1,u_2), (u_2,u_3), \dots, (u_{\ell-1},u_\ell), (u_\ell,u_1).\]
Note that if $G$ has no invertible pair then, for any implication class $I$ of $G$,
$I \cap I^{-1} = \emptyset$. Also note that having an invertible pair is equivalent to
having a circuit of length two, and hence if $G$ has no invertible pair then, for 
any implication class $I$ of $G$ neither $I$ nor $I^{-1}$ contains a circuit of 
length two.

\begin{tm} \label{pair_ordering}
If a reflexive graph $G$ does not contain an invertible pair, then it is a strong
cocomparability graph.
\end{tm}

Before proving the theorem, we give an overview of the proof. 
To show that $G$ is a strong cocomparability graph we prove that it has a strong
cocomparability ordering $\prec$ by constructing the set $U_{\prec}$. The set 
$U_{\prec}$ is constructed by adding for each implication class $I$ either $I$ or
$I^{-1}$. A necessary property for $U_{\prec}$ to have is that it does not 
contain a {\em circuit}. We will show that it is always possible to add either $I$ 
or $I^{-1}$ for each implication class $I$ without creating a circuit. 
For the sake of contradiction, suppose that circuits are created for the first time 
when adding an implication class $I$ or its inverse $I^{-1}$. By considering 
a shortest circuit and two cases depending on the implication classes of the pairs 
in the circuit we will show that the circuit has either one or two pairs from the 
newly added implication class. This will allow us to deduce that a circuit already 
exists even before adding either of $I$ and $I^{-1}$, which is a contradiction.
       
\pf Suppose that $G$ contains no invertible pair. In view of Proposition \ref{scc},
it suffices to show that there exists a set $U$ of pairs of distinct vertices of $G$
which satisfies the following two properties:
\begin{itemize}
\item $U$ contains either $I$ or $I^{-1}$ for all implication classes $I$;
\item $U$ contains no circuit.
\end{itemize}

In fact, we will show such a set $U$ can be constructed by adding implication 
classes one at a time, in no particular order except priority will be given to 
non-trivial implication classes. Initially, $U = \emptyset$. At each step, 
we consider an implication $I$ and its inverse class $I^{-1}$, no pairs of which 
have been added to $U$. We will add all pairs of $I$ to $U$ (that is, replace $U$ 
by $U \cup I$) if $U \cup I$ contains no circuit, otherwise we add all pairs of 
$I^{-1}$ to $U$. It remains to show that if $U$ does not contain a circuit, then 
at least one of $U \cup I$ and $U \cup I^{-1}$ does not contain a circuit. 

Suppose that $U \cup I$ contains circuits. Let $C$ be a shortest circuit contained 
in $U \cup I$. Suppose that $I$ is trivial, i.e., $I = \{(u,v)\}$ for some pair 
$(u,v)$. Since there is no circuit in $U$, $(u,v)$ is in $C$. Since $C$ is 
a shortest circuit in $U \cup I$, $(u,v)$ occurs exactly once in $C$. Similarly, 
if $U \cup I^{-1}$ also contains a circuit, then it has a circuit $C'$ in which
$(v,u)$ occurs exactly once in $C'$. Hence $(C \cup C') \setminus \{(u,v),(v,u)\}$
is a circuit in $U$, contradicting the assumption that $U$ has no circuit. 
Hence $I$ is non-trivial. This implies that, prior to considering $I$ and $I^{-1}$, 
no pair in $U$ forms an trivial implication class. It follows that for 
each $(u,v) \in U \cup I$ there exists a pair $(u',v') \in U \cup I$ such that 
$(u',v') \neq (u,v)$ and $(u',v') \Lambda (u,v)$; in particular, $u'$ is adjacent 
to $u$ but not to $v$ and $v'$ is adjacent to $v$ but not to $u$.

Denote $C:\ (x_1,x_2), (x_2,x_3), \dots, (x_{a-1},x_a), (x_a,x_1)$. Note that 
$a \geq 3$ because $G$ contains no invertible pair. Since $C$ is a shortest circuit
in $U \cup I$, for all $i, j \in \{1, 2, \dots, a\}$, $(x_i,x_j) \in U \cup I$ if 
and only if $j = i+1$ (here and in the rest of proof the subscripts of $x$'s are 
taken modulo $a$). This implies in particular that the pairs on $C$ are pairwise 
distinct. We know from the above that for each $i$ there is a vertex 
adjacent to $x_i$ but not to $x_{i-1}$ and there is also a vertex adjacent to $x_i$
but not to $x_{i+1}$. Below we distinguish two cases depending on whether such 
vertices can be chosen to be the same one for some $i$.

{\bf Case 1.} For some $i$, there is a vertex adjacent to $x_i$ but not
to either of $x_{i-1}$ and $x_{i+1}$. By symmetry we may assume that $i = 1$ and 
$x'_1$ is such a vertex (i.e., adjacent to $x_1$ but not to $x_a$ or $x_2$).  

{\bf Claim 1.1.} For each $i = 2, \dots, a$, let $x'_i$ be a vertex adjacent to 
$x_i$ but not to $x_{i-1}$. Let $z'$ be a vertex adjacent to $x_a$ but not to $x_1$.
Then, for each $i = 1, 2, \dots, a$, $x'_i$ is not adjacent to any of 
$x_a, x_1, \dots, x_{i-1}$ and $x_i$ is not adjacent to any of 
$z', x'_1, x'_2, \dots, x'_{i-1}$ for each $i = 1, 2, \dots, a$,
except when $i = a$ in which case $x_a$ is adjacent both $x'_a$ and $z'$ by 
definition. Moreover, 
\[C':\ (x'_1,x'_2), (x'_2,x'_3), \dots, (x'_{a-1},x'_a), (x'_a,x'_1)\]
is a shortest circuit in $U \cup I$.

{\bf Proof of Claim 1.1.} The statement is clearly true for $i = 1$. By definition
$x'_2$ is not adjacent to $x_1$ and by assumption $x_2$ is not adjacent to $x'_1$. 
If $x'_2$ is adjacent to $x_a$ then $(x_1,x_2) \Lambda (x'_1,x'_2) \Lambda 
(x_1,x_a)$, which implies $(x_1,x_a) \in U \cup I$, a contradiction to the 
observation above. So $x'_2$ is not adjacent to $x_a$. If $x_2$ is adjacent to $z'$
then $(x_1,x_2) \Lambda (x'_1,z') \Lambda (x_1,x_a)$, a contradiction to the
observation above. So $x_2$ is not adjacent to $z'$. Hence the statement is 
also true for $i = 2$. Assume that $2 \leq k < a$ and the statement is true for 
$i = 1, 2, \dots, k$. Suppose that $x_{k+1}$ is adjacent to $x'_k$.  
The vertex $x_{k+1}$ must be adjacent to $x'_{k-1}$ as otherwise 
$(x_{k-1},x_k) \Lambda (x'_{k-1},x'_k) \Lambda (x_{k-1},x_{k+1})$, which 
implies $(x_{k-1},x_{k+1}) \in U \cup I$, a contradiction. Similarly, $x_{k+1}$ 
must be adjacent to $x'_{k-2}$ as otherwise $(x_{k-2},x_{k-1}) \Lambda 
(x'_{k-2},x'_{k-1}) \Lambda (x_{k-2},x_{k+1})$, which implies 
$(x_{k-2},x_{k+1}) \in U \cup I$, a contradiction. Continuing this way we have that
$x_{k+1}$ is adjacent to all of $z', x'_1, x'_2, \dots, x'_k$. 
Let $x'$ be a vertex adjacent to $x_k$ but not to $x_{k+1}$. We know from above 
such a vertex $x'$ exists. Then $x'$ must be adjacent to $x_{k-1}$ as 
otherwise $(x_{k-1},x_k) \Lambda (x'_{k-1},x') \Lambda (x_{k+1},x_k)$, which 
implies $(x_{k+1},x_k) \in U \cup I$, a contradiction. Similarly, $x'$ must be
adjacent to $x_{k-2}$ as otherwise $(x_{k-2},x_{k-1}) \Lambda (x'_{k-2},x') 
\Lambda (x_{k+1},x_{k-1})$, which implies $(x_{k+1},x_{k-1}) \in U \cup I$,
a contradiction. Continuing this way we have $x'$ is adjacent to all of 
$x_a, x_1, x_2, \dots, x_k$. But then $(x_k,x_{k+1}) \Lambda (x',z') \Lambda
(x_1,x_{k+1})$, which implies $(x_1,x_{k+1}) \in U \cup I$, a contradiction.
Hence $x_{k+1}$ is not adjacent to $x'_k$. It cannot be adjacent to 
$x'_{k-1}$ as otherwise $(x_k,x_{k-1}) \Lambda (x'_k,x'_{k-1}) \Lambda 
(x_k,x_{k-1})$, which implies $(x_k,x_{k-1}) \in U \cup I$, a contradiction.
Similarly, it cannot be adjacent to $x'_{k-2}$ as otherwise $(x_k,x_{k+1}) 
\Lambda (x'_k,x'_{k-2}) \Lambda (x_k,x_{k-2})$, which implies $(x_k,x_{k-2}) 
\in U \cup I$, a contradiction. Continuing this way we have that $x_{k+1}$ is 
not adjacent to any of $z', x'_1, x'_2, \dots, x'_k$, except when 
$k+1=a$ in which case $x_a$ is by definition adjacent to $z'$. 

If $x'_{k+1}$ is adjacent to $x_{k-1}$, then $(x_{k-1},x_k) \Lambda (x'_{k+1},x'_k)
\Lambda (x_{k+1},x_k)$, implying $(x_{k+1},x_k) \in U \cup I$, 
a contradiction. So $x'_{k+1}$ is not adjacent to $x_{k-1}$.
If $x'_{k+1}$ is adjacent to $x_{k-2}$, then $(x_{k-2},x_{k-1}) \Lambda 
(x'_{k+1},x'_{k-1}) \Lambda (x_{k+1},x_{k-1})$, which implies 
$(x_{k+1},x_{k-1}) \in U \cup I$, a contradiction. So $x'_{k+1}$ is not adjacent to
$x_{k-2}$. Continuing this way we see that $x'_{k+1}$ is not adjacent to any of 
$x_a, x_1, \dots, x_k$, except when $k+1 = a$ in which case $x'_a$ is by definition
adjacent to $x_a$. 

Since $x_i$ is adjacent to $x'_i$ but not to $x'_{i+1}$ and $x_{i+1}$ is adjacent to
$x'_{i+1}$ but not to $x'_i$, $(x_i,x_{i+1}) \Lambda (x'_i,x'_{i+1})$ for each
$i = 1, 2, \dots, a$. Therefore $C'$ is a circuit in $U \cup I$. Since $C'$ has 
the same length as $C$, it is a shortest circuit in $U \cup I$. 
\qed

It follows from Claim 1.1 that the edges $x_ix'_i$ form an independent set of 
edges in $G$, that is, $x_i$ is not adjacent to $x'_j$ for all $i \neq j$.

{\bf Claim 1.2.} For any vertex $x'$ that is adjacent to $x_i$ but not to $x_{i+1}$,
$x'$ is not adjacent to any of $x_1, \dots, x_{i-1}, x_{i+1}, \dots, x_a$. 
Moreover, 
\[(x'_1,x'_2), \dots, (x'_{i-2},x'_{i-1}), (x'_{i-1},x'), (x',x_{i+1}), 
\dots, (x'_{a-1},x'_a), (x'_a,x'_1)\]
is a shortest circuit in $U \cup I$.

{\bf Proof of Claim 1.2.} If $x'$ is adjacent to $x_{i-1}$, then 
$(x_i,x_{i+1}) \Lambda (x',x'_{i+1}) \Lambda (x_{i-1},x_{i+1})$,   
which implies $(x_{i-1},x_{i+1}) \in U \cup I$, a contradiction. So $x'$ is 
adjacent to $x_i$ but not to $x_{i-1}$ or $x_{i+1}$. By Claim 1.1, $x'$ is not 
adjacent to any of $x_1, \dots, x_{i-1}, x_{i+1}, \dots, x_a$. Since 
$(x_{i-1},x_i) \Lambda (x'_{i-1},x')$, $(x_i,x_{i+1}) \Lambda (x',x'_{i+1})$ and
$(x_j,x_{j+1}) \Lambda (x'_j,x'_{j+1})$ for all $j \neq i-1$ or $i$, 
$(x'_1,x'_2), \dots, (x'_{i-2},x'_{i-1}), (x'_{i-1},x'), (x',x_{i+1}), \dots, 
(x'_{a-1},x'_a), (x'_a,x'_1)$ is a shortest circuit in $U \cup I$..
\qed

{\bf Claim 1.3.} For any vertex $x$ that is adjacent to $x'_i$ but not to 
$x'_{i+1}$, $x$ is not adjacent to any of 
$x'_1, \dots, x'_{i-1}, x'_{i+1}, \dots, x'_a$. Moreover, 
\[(x_1,x_2), \dots, (x_{i-2},x_{i-1}), (x_{i-1},x), (x,x_{i+1}), 
\dots, (x_{a-1},x_a), (x_a,x_1)\]
is a shortest circuit in $U \cup I$.

{\bf Proof of Claim 1.3.} The proof is the same as for Claim 1.2 when the roles of 
$x_i$ and $x'_i$ are switched for each $i$.
\qed

{\bf Claim 1.4.} The circuit $C$ contains a unique pair from $I$. 

{\bf Proof of Claim 1.4.} Since $U$ contains no circuit, there is at least one pair
on $C$ which is from $I$. Without loss of generality assume $(x_a,x_1) \in I$. 
Suppose that $(x_j,x_{j+1}) \in I$ for some $j \neq a$. Then there is a sequence
\[(u_1,v_1) \Lambda (u_2,v_2) \Lambda \cdots \Lambda (u_k,v_k)\]
of pairs in $I$ such that $(u_1,v_1) = (x_a,x_1)$ and $(u_k,v_k) = (x_j,x_{j+1})$.
Moreover, we assume that $(u_i,v_i) \neq (u_{i+1},v_{i+1})$ for each $i$. 
This guarantees that $u_{i+1}$ is adjacent to $u_i$ but not to $v_i$ and 
$v_{i+1}$ is adjacent to $v_i$ but not to $u_i$.
It follows from Claims 1.3 and 1.4, for each even $i$ with $1 \leq i \leq k$, 
$(u_i,v_i), (v_i,x'_2), (x'_2,x'_3), \dots, (x'_{a-2},x'_{a-1}), (x'_{a-1},u_i)$
is a shortest circuit in $U \cup I$. Moreover, $u_i$ is not adjacent to any of
$x_1, x_2, \dots, x_{a-1}$. In particular, when $k$ is even, $u_k$ is not adjacent
to $x_j = u_k$, which is a contradiction. So $k$ is odd. By Claims 1.3 and 1.4, 
$(u_k,v_k), (v_k,x_2), (x_2,x_3), \dots, (x_{a-2},x_{a-1}), (x_{a-1},u_k)$ is
a shortest circuit in $U \cup I$. This is a contradiction because $u_k = x_j$ and
$(u_k,v_k), (v_k,x_2), (x_2,x_3), \dots, (x_{j-1},x_j)$ is a shorter circuit 
in $U \cup I$. 
\qed

{\bf Claim 1.5.} For each pair $(u,v) \in I$, $U \cup I$ contains a circuit in which
$(u,v)$ is the only pair from $I$. 

{\bf Proof of Claim 1.5.} By Claim 1.4, we may again assume that $(x_a,x_1)$ is 
the only pair in $C$ that is from $I$. Since $(x_a,x_1)$ and $(u,v)$ are both
from $I$, there is a sequence 
\[(u_1,v_1) \Lambda (u_2,v_2) \Lambda \cdots \Lambda (u_k,v_k)\]
of pairs in $I$ such that $(u_1,v_1) = (x_a,x_1)$, $(u_k,v_k) = (u,v)$, and
$(u_i,v_i) \neq (u_{i+1},v_{i+1})$ for each $i$.  
The proof of Claim 1.4 also shows that if $k$ is even then substituting $(u,v)$
for the pair $(x'_a,x'_1)$ in $C'$ we obtain a circuit in which $(u,v)$ is 
the only pair from $I$. Similarly, if $k$ is odd then substituting $(u,v)$ for 
the pair $(x_a,x_1)$ in $C$ we obtain a circuit in which $(u,v)$ is the only pair
from $I$.
\qed 

{\bf Case 2.} For each $i$, any vertex adjacent to $x_i$ must be adjacent to at 
least one of $x_{i-1}$ and $x_{i+1}$.

{\bf Claim 2.1.} For each $i = 1, 2, \dots, a$, let $x'_i$ be a vertex adjacent to
$x_i$ but not to $x_{i-1}$ and $x''_i$ be a vertex adjacent to $x_i$ but not to
$x_{i+1}$. Then 
\begin{itemize}
\item $x'_i$ is adjacent to all of $x_a, x_1, x_2, \dots, x_{i-2}$, except 
      when $i = 1$ in which case it is not adjacent to $x_a$.
\item $x''_i$ is adjacent to all of $x_a, x_1, \dots, x_i$, 
      except when $i = a-1$ in which case it is not adjacent to $x_a$ and
      when $i = a$ in which case it is not adjacent to $x_1$.
\item $x_i$ is adjacent to all of $x'_1,x'_2,\dots,x'_{i-1},x''_a,x''_1,x''_2, 
      \dots,x''_{i-2}$, except when $i = 1$ in which case it is not adjacent to 
      $x'_2$ or $x''_a$.
\end{itemize}

{\bf Proof of Claim 2.1.} By definition, $x'_1$ is not adjacent to $x_a$ and $x_1$ 
is not adjacent to $x'_2$ or $x''_a$. The vertex $x''_1$ is adjacent to $x_1$ but 
not to $x_2$ so it is adjacent to $x_a$ by assumption. Thus the 
statement holds for $i = 1$. The vertex $x'_2$ must be adjacent to $x_a$ as 
otherwise $(x_1,x_2) \Lambda (x''_1,x'_2) \Lambda (x_a,x_2)$, which implies 
$(x_a,x_2) \in U\cup I$, a contradiction to the choice of $C$. By assumption $x_2$ 
is adjacent to $x'_1$. If $x_2$ is not adjacent to $x''_a$ then 
$(x_a,x_1) \Lambda (x''_a,x'_1) \Lambda (x_a,x_2)$, a contradiction.
So $x_2$ is adjacent to both $x'_1$ and $x''_a$. Note that $x'_2$ is adjacent to 
$x_2$ but not to $x_1$ so it must be adjacent to $x_3$ by assumption. 
Now consider $x''_2$. It is adjacent to $x_2$ by definition. It must be adjacent to
$x_a$ as otherwise $(x_a,x_1) \Lambda (x'_2,x''_2) \Lambda (x_3,x_1)$, 
a contradiction. Hence the statement holds also for $i=2$. 
Assume that $k \geq 2$ and the statement holds for each $i = 1, 2, \dots, k$. 
We show that it holds for $i = k+1$. 
Consider first the vertex $x'_{k+1}$. It must be adjacent to $x_{k-1}$ as 
otherwise $(x_k,x_{k+1}) \Lambda (x''_k,x'_{k+1}) \Lambda (x_{k-1},x_{k+1})$, 
a contradiction. It must also be adjacent to $x_{k-2}$ as otherwise 
$(x_k,x_{k+1}) \Lambda (x''_k,x'_{k+1}) \Lambda (x_{k-2},x_{k+1})$, 
a contradiction. Continuing this way, we see that $x'_{k+1}$ is adjacent to all of 
$x_a, x_1, x_2, \dots, x_{k-1}$. Now consider $x_{k+1}$. It must be
adjacent to $x'_1$ as otherwise $(x_k,x_{k+1}) \Lambda (x'_1,x'_{k+1}) \Lambda
(x_k,x_a)$, a contradiction. It must also be adjacent to $x'_2$ as otherwise 
$(x_k,x_{k+1}) \Lambda (x'_2,x'_{k+1}) \Lambda (x_k,x_1)$, a contradiction.
Continuing this way we see that $x_{k+1}$ is adjacent to all of 
$x'_1, x'_2, \dots, x'_k$. A similar proof shows that $x_{k+1}$ is adjacent to all 
of $x''_a,x''_1,x''_2, \dots,x''_{k-1}$. Note that $x'_{k+1}$ is adjacent to 
$x_{k+1}$ but not to $x_k$ so it must be adjacent to $x_{k+2}$ by assumption. 
The vertex $x_{k+2}$ must also be adjacent to $x'_k$ as otherwise
$(x_{k-1},x_k) \Lambda (x'_{k+1},x'_k) \Lambda (x_{k+2},x_k)$, a contradiction. 
Continuing this way we see that $x_{k+2}$ is adjacent to all of 
$x'_1, x'_2, \dots, x'_{k+1}$. Finally we consider $x''_{k+1}$. It is adjacent to 
$x_{k+1}$ by definition. It must be adjacent to $x_{k-1}$ as otherwise 
$(x_{k-1},x_k) \Lambda (x'_{k+1}, x''_{k+1}) \Lambda (x_k,x_{k+2})$, 
a contradiction. It must be adjacent to $x_{k-2}$ as otherwise 
$(x_{k-2},x_{k-1}) \Lambda (x'_k,x''_{k+1}) \Lambda (x_{k-1},x_{k+2})$. 
Continuing this manner we see that $x''_{k+1}$ is adjacent to 
all of $x_a, x_1, \dots, x_{k+1}$.
\qed

Note that for each $i = 1, 2, \dots, a$, $x_ix'_{i+1}$ and $x_ix''_{i-1}$ are 
edges in the complement $\overline{G}$ of $G$. It follows from Claim 2.1 that 
in $\overline{G}$, the edges $x_ix'_{i+1}$ form an independent set of edges,
that is, $x_i$ is not adjacent to $x'_j$ for all $i \neq j-1$; similarly,
the edges $x_ix''_{i-1}$ form an independent set of edges, that is, 
$x_i$ is not adjacent to $x''_j$ for all $i \neq j+1$.
 
{\bf Claim 2.2.} For each $i = 1, 2, \dots, a$, let $x'_i$ be a vertex adjacent to
$x_i$ but not to $x_{i-1}$ and $x''_i$ be a vertex adjacent to $x_i$ but not to
$x_{i+1}$. Then 
\[C'':\ (x'_a,x'_{a-1}), (x'_{a-1},x'_{a-2}), \dots, (x'_2,x'_1), (x'_1,x'_a)\]
is a shortest circuit in $U \cup I$.  

{\bf Proof of Claim 2.2.} By Claim 2.1, each $x_i$ is adjacent to $x'_{i+2}$ and 
$x_{i+1}$ is adjacent to $x'_{i+1}$. So $(x_i,x_{i+1}) \Lambda (x'_{i+2},x'_{i+1})$
for each $i$. Hence $C''$ is a circuit in $U \cup I$. It is the shortest because it
has the same length as $C$.
\qed 

{\bf Claim 2.3.} Let $i$ be fixed. If $(x_i,x_{i+1}) \Lambda (u,v)$ for some 
$(u,v) \neq (x_i,x_{i+1})$, then substituting $(u,v)$ for $(x'_{i+2},x'_{i+1})$ in 
$C''$ results in a shortest circuit in $U \cup I$. Similarly, if 
$(x'_{i+2},x'_{i+1}) \Lambda (u',v')$ for some $(u',v') \neq (x'_{i+2},x'_{i+1})$, 
then substituting $(u',v')$ for $(x_i,x_{i+1})$ in $C$ results in a circuit in 
$U \cup I$.

{\bf Proof of Claim 2.3.} If $(x_i,x_{i+1}) \Lambda (u,v)$, then $u$ can play 
the role of $x'_{i+2}$ and $v$ can play role of $x'_{i+1}$ in Claim 2.1. 
Hence by Claim 2.2 substituting $(u,v)$ for $(x'_{i+2},x'_{i+1})$ in $C''$ results 
in a circuit in $U \cup I$. The new circuit is a shortest one in $U \cup I$ because
it has the same length as $C$. The proof of the second part is the same when 
the roles of $C$ and $C''$ are switched in both Claims 2.1 and 2.2. 
\qed
 
{\bf Claim 2.4.} $C$ contains one or two pairs from $I$. 

{\bf Proof of Claim 2.4.} Since there is no circuit in $U$, $C$ must contain at
least one pair from $I$. Suppose that $(x_r,x_{r+1}), (x_s,x_{s+1}) \in I$ for some
$r \neq s$. Since both $(x_r,x_{r+1}), (x_s,x_{s+1})$ are from $I$, there is 
a sequence  
\[(u_1,v_1) \Lambda (u_2,v_2) \Lambda \cdots \Lambda (u_k,v_k)\]
of pairs in $I$ such that $(u_1,v_1) = (x_r,x_{r+1})$, $(u_k,v_k) = (x_s,x_{s+1})$,
and $(u_i,v_i) \neq (u_{i+1},v_{i+1})$ for each $i$. If $k$ is odd, then it follows
from Claim 2.3 substituting $(x_s,x_{s+1})$ for $(x_r,x_{r+1})$ in $C$ results in 
a shortest circuit in $U \cup I$. But the new circuit has equal vertices $x_s$ (and
$x_{s+1}$) so it is not the shortest, a contradiction. Hence $k$ must be even.
Therefore any two pairs in $C$ which are from $I$ can only be connected by such 
a sequence consisting of an even number of pairs. Suppose that $C$ contains three
pairs from $I$. Then there is a sequence as above consisting of an even number
of pairs connecting the first and the second pairs, and there is a sequence 
consisting consisting of an even number of pairs connecting the second and the third
pairs. By concatenating the two sequences we obtain a sequence consisting of an
odd number of pairs connecting the first and the third pairs, which contradicts
the above fact.
\qed

{\bf Claim 2.5.} If $C$ contains exactly one pair from $I$, then for any 
$(u,v) \in I$, $U \cup I$ contains a circuit in which $(u,v)$ is the only pair from
$I$. If $C$ contains exactly two pairs from $I$, then for any two pairs 
$(u_1,v_1), (u_k,v_k) \in I$ connected by a sequence
\[(u_1,v_1) \Lambda (u_2,v_2) \Lambda \cdots \Lambda (u_k,v_k)\]
where $(u_i,v_i) \neq (u_{i+1},v_{i+1})$ and $k$ is odd, $U \cup I$ contains
a circuit in which $(u_1,v_1), (u_k,v_k)$ are the only pairs from $I$.

{\bf Proof of Claim 2.5.} If $C$ contains exactly one pair from $I$, then $C''$
also contains exactly one pair from $I$. By Claim 2.3, substituting $(u,v)$ for
the pair in $C$ or the pair in $C''$ results in a circuit in which $(u,v)$ is 
the only pair from $I$. Similarly, if $C$ contains exactly two pairs from $I$,
then $C''$ also contains exactly two pairs from $I$. By Claim 2.3, substituting
$(u_1,v_1), (u_k,v_k)$ for the two pairs in $C$ or the two pairs in $C''$ results
in a circuit in which $(u_1,v_1), (u_k,v_k)$ are the only two pairs from $I$.  
\qed

We know from above that if $U \cup I$ contains a circuit then Case 1 or Case 2 
occurs for any shortest circuit in $U \cup I$. In any case, either $U \cup I$ 
has a circuit which contains exactly one pair from $I$ or a circuit which contains
exactly two pairs from $I$. The same is true for $U \cup I^{-1}$.
Suppose that both $U \cup I$ and $U \cup I^{-1}$ contain circuits. 
We show that $U$ has circuits, which contradicts our assumption.
Let $(u,v), (u',v') \in I$ be connected by the sequence 
\[(u,v) = (u_1,v_1) \Lambda (u_2,v_2) \Lambda \cdots \Lambda (u_k,v_k) = (u',v')\]
where $(u_i,v_i) \neq (u_{i+1},v_{i+1})$ and $k$ is odd. 
When each of $U \cup I$ and $U \cup I^{-1}$ has a circuit which contains exactly
one pair from $I$, by Claim 1.5 $U \cup I$ has a circuit $C$ in which 
$(u,v)$ is the only pair from $I$ and $U \cup I^{-1}$ has a circuit $C'$ in which 
$(v,u)$ is the only pair from $I^{-1}$. Then $(C \cup C') \setminus \{(u,v),(v,u)\}$
is a circuit in $U$.
When each of $U \cup I$ and $U \cup I^{-1}$ has a circuit which contains 
exactly two pairs from $I$, by Claim 2.5 $U \cup I$ has a circuit $C$ in which  
$(u,v), (u',v')$ are the only pairs from $I$ and $U \cup I^{-1}$ has a circuit $C'$ 
in which $(v,u), (v',u')$ are the only pairs from $I^{-1}$. Then 
$(C \cup C') \setminus \{(u,v),(u',v'),(v,u),(v',u')\}$ is a union of two circuits 
in $U$.
Suppose that one of $U \cup I$ and $U \cup I^{-1}$ has a circuit which contains 
exactly one pair from $I$ and the other has a circuit $C''$ which contains exactly 
two pairs from $I$. Without loss of generality assume that $U \cup I$ has a circuit
which contains exactly one pair from $I$. By Claims 1.5 and 2.5, 
$U \cup I$ has a circuit $C$ in which $(u,v)$ is the only pair from $I$ and 
a circuit $C'$ in which $(u',v')$ is the only pair from $I$, and $U \cup I^{-1}$
has a circuit $C''$ in which $(v,u), (v',u')$ are the only pairs from $I$. 
Then $(C \cup C' \cup C'') \setminus \{(u,v), (u',v'), (v,u), (v',u')\}$
is a circuit in $U$. Therefore at least one of $U \cup I$ and $U \cup I^{-1}$ 
does not contain a circuit. This completes the proof.
\qed

Combining Propositions \ref{ip}, \ref{aux}, and Theorem \ref{pair_ordering}, we 
have the following: 

\begin{tm} \label{main}
The following statements are equivalent for a reflexive graph $G$:
\begin{description}
\item{(i)} $G$ is a strong cocomparability graph;
\item{(ii)} $G$ contains no invertible pairs;
\item{(iii)} $G^+$ is bipartite.
\qed
\end{description}
\end{tm}

Given a reflexive graph $G$, in polynomial time one can construct $G^+$ and test
whether it is bipartite. Thus we have the following:

\begin{cor} \label{ip-test}
There is a polynomial-time algorithm to decide whether a reflexive graph is 
a strong cocomparability graph.
\qed
\end{cor}

\section{Concluding remarks and open problems}

Strong cocomparability graphs form a proper subclass of cocomparability graphs.
Gilmore and Hoffman \cite{gh} proved that a reflexive graph is an interval graph 
if and only if it is both a chordal graph and a cocomparability graph. 
The following may be viewed as a strengthening of the theorem of 
Gilmore and Hoffman.

\begin{tm} 
A reflexive graph $G$ is an interval graph if and only if it is both a strongly 
chordal graph and a strong cocomparability graph.
\end{tm}

\pf Suppose that $G$ is an interval graph. Then by Theorems \ref{gamma} and
\ref{gamma-slash} it is strongly chordal. Since $G$ is an interval graph, it has 
an interval ordering, which is a vertex ordering $\prec$ such that for all 
$x \prec y \prec z$, $xz \in E(G)$ implies $xy \in E(G)$. It is easy to check that 
$\prec$ is a strong cocomparability ordering. Hence $G$ is also a strong 
cocomparability graph. 

On the other hand, suppose that $G$ is not an interval graph. Then according to
\cite{lb} it contains an induced cycle of length $\geq 4$ or an asteroidal triple
(i.e., a set of three vertices such that any two of them are joined by a path which
contains no neighbour of the third vertex).
If it contains an induced cycle of length $\geq 4$ then it is not strongly chordal.
Suppose that $u, v, w$ form an asteroidal triple in $G$. Then there is a path $P$
joining $u, v$ which contains no neighbour of $w$. For each edge $xy$ on $P$, 
$xwy$ is an induced $P_3$ in $\overline{G}$. By Proposition \ref{pc}, 
$(x,w) \Lambda (y,w)$. It follows that $(u,w) \sim (v,w)$. Similarly, 
$(u,v) \sim (w,v)$ and $(w,u) \sim (v,u)$. Hence 
\[(u,v) \sim (w,v) \sim (w,u) \sim (v,u)\]
which shows that $u, v$ are an invertible pair in $G$ and, by Theorem \ref{main}, 
$G$ is not a strong cocomparability graph. 
\qed

A forbidden induced subgraph characterization of cocomparability graphs has 
been given by Gallai \cite{gallai}. It would be very interesting to find a forbidden 
induced subgraph characterization for strong cocomparability graphs. 

\begin{prob}
Characterize strong cocomparability graphs by forbidden subgraphs.
\end{prob}

Gallai's forbidden induced subgraph characterization of cocomparability graphs
has 7 infinite families (including for example cycles with at least five vertices) 
and 10 of other subgraphs, each on seven vertices. In addition to these, we have 
found by a computer search that there are many other minimal graphs that are not 
strong cocomparability graphs (but are cocomparability graphs) -- 5 on six vertices, 
5 on seven vertices, 12 on eight vertices, and 30 on nine vertices.

Thus the above problem seems difficult at the moment. On the other hand, there
is a simple forbidden subgraph characterization for strong cocomparability graphs
that are bipartite. 

A reflexive graph is {\em bipartite} if the graph obtained from it by deleting all loops is 
bipartite in the traditional sense. We use other terminology on bigraphs and bipartite 
graphs in the same sense, by applying them to the graph with loops removed; notably
this is how we define the class of chordal bigraphs for the reflexive case.

Figure 2 depicts the forbidden subgraphs $S$ and $T$ for bipartite strong cocomparability 
graphs.

\begin{figure}[ht]
\begin{center}
\begin{tikzpicture}[>=latex]
\node [style=blackvertex] (1) at (-6,0) {};
                \node [style=blackvertex] (2) at (-4,0) {};
                \node [style=blackvertex] (3) at (-2,0) {};
                \node [style=blackvertex] (4) at (-2,1.5) {};
                \node [style=blackvertex] (5) at (-4,1.5) {};
                \node [style=blackvertex] (6) at (-6,1.5) {};
                \draw [ ] (1)--(2)--(3)--(4)--(5)--(6);
                \draw [ ] (2) to (5);
                
\node [style=blackvertex] (1) at (0,0) {};
                \node [style=blackvertex] (2) at (2,0) {};
                \node [style=blackvertex] (3) at (4,0) {};
                \node [style=blackvertex] (4) at (4,1.5) {};
                \node [style=blackvertex] (5) at (2,1.5) {};
                \node [style=blackvertex] (6) at (0,1.5) {};
                \node [style=blackvertex] (7) at (6,1.5) {};
                \draw [ ] (6)--(5)--(2)--(3);
                \draw [ ] (1)--(6);
                \draw [ ] (5)--(4)--(7);
\end{tikzpicture}
\end{center}
\vspace{-2mm}
\caption{The two minimal forbidden subgraphs, $S$ on the left and $T$ on the right.}
\label{6}
\end{figure}
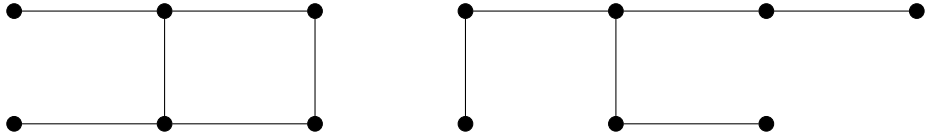

\begin{tm}\label{bipartite}
A bipartite graph is a strong cocomparability graph if and only is it is a chordal 
bigraph that does not contain $S$ or $T$ as a subgraph.
\end{tm}

Note that the characterization is stated for forbidden (not necessarily induced) subgraphs. 
It can of course be translated to a forbidden induced subgraph characterization, by forbidding 
all bipartite supergraphs of $S$ and $T$. 

A {\em caterpillar} is a tree $G$ such that after the deletion of all its leaves (vertices of degree one) 
$G$ becomes a path, called the {\em spine} of $G$. If a vertex $v$ on the spine of a caterpillar $G$
has degree two in $G$, we say we {\em replicate} it if we add to $G$ additional vertices $v', v'', \dots, v^*$ 
of degree two, with exactly the same neighbours as $v$; we call $v$ the {\em shadow} of all the
replicated vertices. Consider a complete bipartite subgraph $K_{2,s}$ (for any $s \ge 2$) of a graph 
$G$, with the bipartition $A \cup B$, where $|A|=2, |B|=s \ge 2$. If all $s$ vertices in $B$ have degree 
two in $G$, we call the subgraph a {\em square in $G$} and the two vertices in $A$ its {\em connectors}.

\pf 
It follows from \cite{gallai} that $T$ or any cycle with at least six vertices is not a cocomparability 
(and hence not a strong cocomparability) graph. Moreover, any spanning bipartite supergraph of 
$T$ contains a bipartite supergraph of $S$. Finally, all bipartite supergraphs of $S$ fail to be strong 
cocomparability graphs; this can be seen by a simple but tedious argument. Therefore, each bipartite 
strong cocomparability graph is a chordal bigraph that does not contain $S$ or $T$ as a subgraph.

To prove the converse, let $G$ be a connected chordal bigraph that does not contain $S$ or $T$ 
as a subgraph. If $G$ is a tree, then the absence of $T$ implies that $G$ is a caterpillar, with spine 
$v_1, v_2, \dots, v_k$. We construct a strong cocomparability ordering of any caterpillar $G$ as 
$v_1 \prec S_1 \prec v_2 \prec S_2 \prec \dots \prec v_k \prec S_k$,  where $S_i$ is the set of 
leaves adjacent to $v_i$, listed in arbitrary order.

If $G$ is not a tree, then it must contain a four-cycle. Moreover, each four-cycle must have two opposite 
vertices of degree two, else we have a copy of $S$. Suppose $abcda$ is a four-cycle with $a, c$ of degree 
two. Then any vertex $x$ adjacent to both $b$ and $d$ must also be of degree two, and by adjoining all 
such vertices $x$ we obtain a square in $G$. This applies even in the case the square is just a four-cycle, 
as in the above example, where $b, d$ are the connectors. (Note that in such a case, the connectors can 
also be of degree two.) Two squares can meet in at most one vertex, which must be a connector in both 
squares, otherwise $G$ had a copy of $S$. Now consider the graph $G^*$ obtained from $G$ by deleting 
all but one of the degree two vertices that are not connectors in any square. We can view the removed 
vertices as being the results of replication, and hence the remaining vertex is their shadow. Since the 
graph is connected and chordal, $G^*$ is a tree, and hence a caterpillar. Moreover, each shadow vertex 
is on the spine, and has no children. The strong cocomparability ordering of the caterpillar $G^*$, as 
described above, can now be modified to a strong cocomparability ordering of $G$ by inserting each of 
the replicated vertices next to their shadow vertex from its square. Because the shadow vertices have 
no children off the spine, this is indeed a strong cocomparability ordering.
\qed

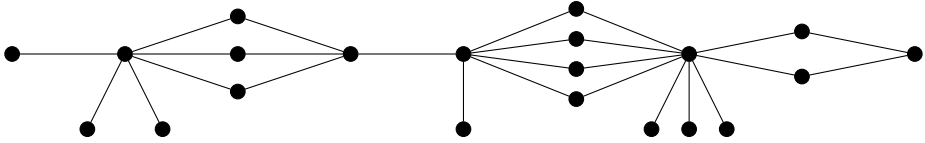
\begin{figure}[ht]
\begin{center}
\begin{tikzpicture}[>=latex]
\node [style=blackvertex] (1) at (-6,0) {};
                \node [style=blackvertex] (2) at (-4.5,0) {};
                \node [style=blackvertex] (3) at (-5,-1) {};
                \node [style=blackvertex] (4) at (-4,-1) {};
                \node [style=blackvertex] (5) at (-3,.5) {};
                \node [style=blackvertex] (6) at (-3,0) {};
                \node [style=blackvertex] (7) at (-3,-.5) {};
                \node [style=blackvertex] (8) at (-1.5,0) {};
                \node [style=blackvertex] (9) at (0,0) {};
                \node [style=blackvertex] (10) at (0,-1) {};
                \node [style=blackvertex] (11) at (1.5,.6) {};
                \node [style=blackvertex] (12) at (1.5,.2) {};
                \node [style=blackvertex] (13) at (1.5,-.2) {};
                \node [style=blackvertex] (14) at (1.5,-.6) {};
                \node [style=blackvertex] (15) at (3,0) {};
                \node [style=blackvertex] (16) at (2.5,-1) {};
                \node [style=blackvertex] (17) at (3,-1) {};
                \node [style=blackvertex] (18) at (3.5,-1) {};
                \node [style=blackvertex] (19) at (4.5,.3) {};
                \node [style=blackvertex] (20) at (4.5,-.3) {};
                \node [style=blackvertex] (21) at (6,0) {};

                \draw [ ] (1)--(2)--(3);
                \draw [ ] (4)--(2)--(5)--(8)--(6)--(2)--(7)--(8)--(9)--(10);
                \draw [ ] (16)--(15)--(11)--(9)--(12)--(15)--(13)--(9)--(14)--(15)--(17);
                \draw [ ] (18)--(15)--(19)--(21)--(20)--(15);
\end{tikzpicture}
\end{center}
\vspace{-2mm}
\caption{A typical example of a bipartite strong cocomparability graph.}
\label{bi-caterpillar}
\end{figure}

We observe that the proof gives an implicit description of all bipartite strong cocomparability graphs -- 
they are obtained from caterpillars by replicating an arbitrary set of non-consecutive degree two vertices.
This is illustrated in Figure \ref{bi-caterpillar}.

Since bipartite strong cocomparability graphs do not contain asteroidal triples, 
they are bipartite permutation graphs \cite{sbs}. On the other hand if a bipartite 
permutation graph does not contain $S$ as a subgraph then it does not contain $T$ 
as a subgraph and hence by Theorem \ref{bipartite} it is a strong cocomparability
graph. 

\begin{cor} 
A bipartite graph is a strong cocomparability graph if and only if it is a bipartite
permutation graph which does not contain $S$ as a subgraph.
\qed
\end{cor}  

We also want point out a relationship between bipartite strong cocomparability 
graphs and cocomparability bigraphs. Clearly, if a cocomparability bigraph is a 
strong cocomparability graph then it is a bipartite strong cocomparability graph. 
Suppose that $G$ is a bipartite strong cocomparability graph. Then by definition 
the adjacency matrix of $G$ has a symmetric $Slash$-free ordering, which implies 
that the biadjacency matrix of $G$ has a $Slash$-free ordering and hence $G$ is a 
cocomparability bigraph. Therefore bipartite strong cocomparability graphs are
precisely those strong cocomparability graphs which are also cocomparability 
bigraphs. 

On a final note, we observe that by complementation we can introduce the class of 
strong comparability graphs. We define $G$ to be a {\em strong comparability graph},
if its complement $\overline{G}$ is a strong cocomparability graph. 
The complementation includes the loops, thus strong comparability graphs are 
irreflexive, and their adjacency matrices have $0$-diagonal.
A square symmetric $0,1$ matrix with $0$-diagonal has a symmetric $I_2$-free 
ordering if and only if it is the adjacency matrix of a strong comparability graph.
The class of strong comparability graphs has an ordering characterization, 
a forbidden structure characterization, and a polynomial recognition algorithm, 
all via the our corresponding results for their complements, i.e., the strong 
cocomparability graphs.

To further motivate these concepts, we recall from \cite{hhlm} the definition of 
a {\em principal} $\Gamma$ (respectively, $Slash$) submatrix of a symmetric matrix 
with $1$-diagonal: it is a $\Gamma$ (respectively, $Slash$) submatrix with one of 
its $1$'s lying on the main diagonal. Then 

\begin{itemize}
\item
a reflexive graph is chordal if and only if its adjacency matrix admits a symmetric ordering without
a principal $\Gamma$ submatrix; and is strongly chordal if and only if its adjacency matrix admits
a symmetric ordering without any $\Gamma$ submatrix, and
\item
a reflexive graph is a cocomparability graph if and only if its adjacency matrix admits a symmetric
ordering without a principal $Slash$ submatrix; and is a strong cocomparability graph if and only if
its adjacency matrix admits a symmetric ordering without any $Slash$ submatrix.
\end{itemize}

Similarly, we define a {\em principal} $I_2$ submatrix in a symmetric matrix with $0$-diagonal:
this is an $I_2$ submatrix with either of the $0$'s lying on the main diagonal. Then 
\begin{itemize}
\item
an irreflexive graph $G$ is a comparability graph if and only if its adjacency matrix admits a symmetric
ordering without a principal $I_2$ submatrix, and is a strong comparability graph if and only if its adjacency
matrix admits a symmetric ordering without any $I_2$ submatrix.
\end{itemize}

These analogies suggest that these concepts are natural and may lead to interesting questions. As a
first example, consider the domination problem, which is known to be NP-complete for chordal graphs
but polynomial-time solvable for strongly chordal graphs \cite{farber}. The domination problem is also
NP-complete for comparability graphs \cite{john}, so we wonder about strong comparability graphs.

\begin{prob}
Determine the complexity of the domination problem for strong comparability graphs.
\end{prob}

We note that for cocomparability graphs the domination problem is polynomial-time solvable \cite{ks}.


\begin{thebibliography}{100}

\bibitem{af} R.P. Anstee and M. Farber, Characterizations of totally
balanced matrices,
J. Algorithms 5 (1984) 215 - 230.

\bibitem{bd} V.L. Beresnev and A.I. Davydov, On matrices with the connectedness 
property, Upravlyaemye sistemy 19 (1979) 3-13 [in Russian].

\bibitem{signed} J. Bok, R. Brewster, P. Hell, N. Jedli\v ckov\' a, and A. Rafiey,
Min orderings and list homomorphism dichotomies for signed and unsigned graphs,
arXiv:2206.01068 [math.CO].

\bibitem{dama}  P. Damaschke, Forbidden ordered subgraphs, 
Topics in Combinatorics and Graph
Theory (1990) 219 - 229.

\bibitem{far} M. Farber, Characterizations of strongly chordal graphs, 
Discrete Math. 43 (1983) 173 - 189. 

\bibitem{farber} M. Farber,
Domination, independent domination, and duality in strongly chordal graphs,
Discrete Appl. Math. 7 (1984) 115 - 130.

\bibitem{fhh} T. Feder, P. Hell, and J. Huang, List homomorphisms and
circular arc graphs, Combinatorica 19 (1999) 487 - 505.

\bibitem{adjusted} T. Feder, P. Hell, J. Huang, and A. Rafiey, Interval graphs, 
adjusted interval digraphs, and reflexive list homomorphisms, Discrete Appl. Math.
160 (2012) 697 - 707.

\bibitem{cags} M. Francis, P. Hell, and J. Stacho,
Forbidden structure characterization of circular-arc graphs and a certifying recognition algorithm,
SODA 2015.

\bibitem{gallai} T. Gallai, Transitiv orientierbare graphen,
Acta Mathematica Academiae Scientiarum Hungarica 18 (1967) 25 - 66.

\bibitem{gh} P.C. Gilmore and A.J. Hoffman, A characterization of comparability
graphs and interval graphs, 
Canad. J. Math. 16 (1962) 539 - 548.


\bibitem{gg} M.C. Golumbic and C.F. Goss, Perfect elimination and chordal
bipartite graphs. J. Graph Theory 2 (1978) 155 - 163.


\bibitem{hchl} P. Hell, C. Hern\' andez-Cruz, J. Huang, and J. C.-H. Lin,
Strong chordality of graphs with possible loops, SIAM Journal on Discrete Math.
35 (2021) 362 - 375.

\bibitem{hhlm} P. Hell, J. Huang, J.C.-H. Lin, and R.M. McConnell,
Bipartite analogues of comparability and cocomparability graphs, 
SIAM J. Discrete Math. 34 (2020) 1969 - 1983.

\bibitem{hhmr} P. Hell, J. Huang, R.M. McConnell, and A. Rafiey,
Min-orderable digraphs,
SIAM J. Discrete Math. 34 (2020) 1910 - 1724.

\bibitem{mastro} P. Hell, M. Mastrolilli, M.M. Nevisi, and A. Rafiey, 
Approximation of minimum cost homomorphisms, ESA, 2012, pp. 587 - 598.

\bibitem{boj} P. Hell, B. Mohar, and A. Rafiey,
Ordering without forbidden patterns,
ESA 2014, LNCS 8737 (2014) 554 - 565.

\bibitem{mono} P. Hell and A. Rafiey,
Monotone proper interval digraphs and min-max orderings,
SIAM J. Discrete Math. 26 (2012) 1576 - 1596.

\bibitem{hsk} A.J. Hoffman, M. Sakarovich, and A. Kolen,
Totally balanced and greedy matrices,
SIAM J. Algebraic Discrete Methods 6 (1985) 721 - 730.

\bibitem{huang} J. Huang, Representation characterizations of chordal
bipartite graphs, J. Combinatorial Theory B 96 (2006) 673 - 683.

\bibitem{huang-oo} J. Huang, Non-edge orientation and vertex ordering
characterizations of some classes of bigraphs,
Discrete Applied Math. 245 (2018) 190 - 193.

\bibitem{hks} A.J. Hoffman, A.W.. Kolen, and M. Sakarovitch, Totally-balanced
and greedy matrices, 
SIAM J. Algebraic and Discrete Methods 6 (1985) 721 - 730.

\bibitem{john} D.S. Johnson,
The NP-completeness column: an ongoing guide,
J. Algorithms 6 (1985) 145 - 159, 291 - 305, 434 - 451.

\bibitem{krw} B. Klinz, R. Rudolf, and G.J. Woeginger,
Permuting matrices to avoid forbidden submatrices,
Discrete Applied Math. 60 (1995) 223 - 248.

\bibitem{ks} D. Kratsch and L. Stewart,
Domination on cocomparability graphs,
SIAM J. Discrete Math. 6 (1993) 400 - 417.

\bibitem{lb} C.G. Lekkerkerker and J.Ch. Boland, 
Representation of a finite graph by a set of intervals on the real line,
Fund. Math. 51 (1962) 45 - 64.

\bibitem{lubiw1} A. Lubiw, $\Gamma$-free matrices,
Master's thesis, University of Waterloo, 1982.

\bibitem{lub} A. Lubiw, Doubly lexical orderings of matrices,
SIAM J. Comput. 16 (1987) 854 - 879.


\bibitem{pt} R. Paige and R.E. Tarjan, Three partition refinement algorithms, 
SIAM .I. Comput. 16 (1987) 973 - 989.

\bibitem{stu} A.M. Shresta, S. Tayu, and S. Ueno, On two-directional orthogonal 
ray graphs, in Proceedings of the 2010 IEEE International Symposium on Circuits and
Systems, IEEE, pp. 1807 - 1810.

\bibitem{spin} J.P. Spinrad, Efficient Graph Representations, 
Fields Inst. Monogr. 19, AMS, Providence, RI, 2003.

\bibitem{sbs} J.P. Spinrad, A. Brandst\"adt and L. Stewart, 
Bipartite permutation graphs, 
Discrete Appl. Math. 18 (1987) 279 - 292.

\end{thebibliography}
\end{document}